\documentclass[12pt]{amsart}
\usepackage{amssymb,amsmath,amsthm}
\newtheorem{theorem}{Theorem}
\newtheorem{lemma}[theorem]{Lemma}
\newtheorem{corollary}[theorem]{Corollary}

\theoremstyle{remark}

\newtheorem*{remark}{Remark}

\numberwithin{theorem}{section} \numberwithin{equation}{section}

\newcommand{\R}{\mathbb{R}}
\newcommand{\C}{\mathbb{C}}

\newcommand{\Q}{\mathbb{Q}}

\newcommand{\Z}{\mathbb{Z}}
\newcommand{\N}{\mathbb{N}}

\def\H{{\mathfrak H}}

\begin{document}
\title{Miyawaki's $\text{F}_{12}$ Spinor L-function Conjecture}
\author{Bernhard Heim}
\address{Max-Planck Institut  f\"ur Mathematik, Vivatsgasse 7, 53111 Bonn, Germany}
\email{heim@mpim-bonn.mpg.de}
\subjclass[2000] {11F}
%05A17 }
%\dedicatory{Very preliminary}
%%\date{\today}
\begin{abstract}
In this paper we prove the Miyawaki conjecture
related to the spinor $L$--function of a Siegel cusp form of weight $12$
and degree $3$ as a special example of results related to Miyawaki lifts of odd degree.
\end{abstract}
\maketitle
\section{Introduction} 
At a time when no one had a clear idea how to systematically generalize
{\it Saito-Kurokawa liftings} to cuspidal automorphic forms, Isao
Miyawaki \cite{Mi92} made two precise 
conjectures. We follow Miyawaki's notation. Let $F_{12}$ be the unique (up to scalar)
Siegel cusp form of degree $3$ and weight $12$. 
Armed with rich numerical data and certain insights,
Miyawaki determined
the local factors of the spinor $L$--function $L(s, F_{12})$ at the
primes $p=2,3$ (see Theorem 4.2 in his paper). 
He further showed that this implies a certain degeneration of the
standard $L$--function of $F$ (\cite{Mi92}, Section 6). 
The conjecture related to the standard $L$--function was recently
proven by Ikeda \cite{Ik06}. 
Miyawaki states in his paper:
\begin{center}
\ldots the author believes that the above theorem (Theorem 4.2)\\will
be true for any 
prime $p$ (Conjecture 4.3).
\end{center}

In this paper we can manifest his belief in a theorem.
Our proof builts on results of Andrianov \cite{An67}, Yamazaki \cite{Ya86},
Ikeda \cite{Ik01}\cite{Ik06}, Hayashida \cite{Ha07} and others.
\\
\\
{\bf Theorem} {\it Let $F_{12}$ be the Siegel cusp form of
  degree $3$ and weight $12$, as  
given in \cite{Mi92}. Then the spinor $L$--function $L(s,F_{12})$ is given
by 
\begin{equation}
L(s,F_{12}) = L(s-9,\Delta) \, L(s-10,\Delta) \,\, L(s, \Delta \otimes g_{20}).
\end{equation}
Here $L(s,\Delta)$ is the Hecke $L$--function of the Ramanujan function $\Delta$ and
$L(s,\Delta \otimes g_{20})$ the Rankin $L$--function of 
the convolution of $\Delta$ with the newform $g_{20}$ of weight $20$.}
\\
\\
Supplementary we want to mention, that the spinor L-function $L(s,F_{12})$ 
has a holomorphic continuation to the whole complex plane and is the 
first example of a spinor L-function of a Siegel cusp forms of degree $3$
which satisfies a functional equation.
This is compatible with Miyawaki's \cite{Mi92} and Andrianov's \cite{An74} conjectures about the
$\Gamma$-factors and functional equation, where 
$s \mapsto 3 \kappa -5-s$ with weight $\kappa$.
A manuscript related to 
the analytic properties and Deligne's conjectures
\cite{De79}, \cite{Pa94}, \cite{Yo01} of the
spinor L-function of cusp forms of odd degree 
is under preparation.
\section*{Acknowledgements}
\noindent
After this paper was finished Prof. Katsurada kindly informed the author about the work of Kenji Murakawa \cite{Mu02}.
We add some remarks at the end of the paper towards the relation.
%%
%%
%%
%\newpage

\section{Ikeda's lifts and Hayashida's construction}
Let $R$ be any subring of the real numbers $\R$. Then
$$
\text{G}^{+}\text{Sp}_n(R):= \big\{ \gamma \in \text{GL}_{2n}(R) \vert \,\, 
\gamma \, J_n \, \gamma^t = \mu(\gamma)\, J_n, \, \mu(\gamma)>0 \big\},
$$
where $J_n = \left( \begin{smallmatrix} 0 & -1_n \\ 1_n & 0 \end{smallmatrix} \right)$.
The group $\text{G}^{+}\text{Sp}_n(\R)$ acts on the Siegel upper half-plane $\H_n$ of degree $n$
by $
\left( \begin{smallmatrix} a & b \\ c & d \end{smallmatrix} \right)(z):=
\left( az +b\right)  \left(cz +d\right)^{-1}$.
Let $n,k$
be positive integers.
For a holomorphic complex valued function $F$ on
$\H_n$ and $\gamma \in \text{G}^{+}\text{Sp}_n(\R)$
the Petersson slash operator $\vert_k$ is given by
$$
F\vert_k \gamma (z):= \det(cz +d)^{-k} \,\, F(\gamma (z)).$$
The space of Siegel modular forms $M_k^n$ of weight $k$ and degree $n$ is given
by all $F$ with the invariance property $F\vert_k \gamma = F$ for $\gamma \in \Gamma_n = \text{Sp}_n(\Z)$
and certain well known growth conditions in the case $n=1$.
The subspace of cusp forms we denote by $S_k^n$. For $F,G \in S_k^n$ we
have the Petersson scalar product $\langle F,G \rangle$.
Let $\mathcal{H}^n$ be the Hecke algebra assoziated to the Hecke pair
$$
\big(   \Gamma_n, \text{M}_{2n}(\Z) \cap \text{G}^{+}\text{Sp}_n(\Q)\big).
$$
Then it is well known that this Hecke algebra decomposes in an infinite restricted tensor product
$\otimes_p \mathcal{H}_p^n$ of local Hecke algebras generated by
\begin{eqnarray*}
T(p)^{(n)} &:= & \Gamma_n \text{diag}\left( 1_n ; p 1_n \right) \Gamma_n, \\
T_r(p^2)^{(n)} & := & \Gamma_n \text{diag}\left( p1_r,1_{n-r}; p1_r, p^2 1_{n-r}\right) \Gamma_n \quad ( 0 \leq r \leq n ).
\end{eqnarray*}
The action of the Hecke algebra on $S_k^n$ is induced by the action of 
$T= \Gamma_n g \Gamma_n = \cup_i \Gamma_n g_i$ (disjoint) via
\begin{equation}
F\vert_k T := \sum_i \mu(g_i)^{nk-n(n+1)/2} \,\, F\vert_k g_i.
\end{equation}
Let $F\in S_k^n$ be an eigenform for all 
Hecke operators $T$ with eigenvalue $\lambda_F(T)$.
For $T(p^d)^{(n)}:= \sum_{g, \mu(g)=p^d} \Gamma_n g \Gamma_n$,
where $g$ runs through $\text{M}_{2n}(\Z) \cap \text{G}^{+}\text{Sp}_n(\Q)$ 
with different elementary divisors, we set $\lambda_F(T) = \lambda_F(p^d)$.

For the reader's convenience we use Miyawaki's notation for the spinor
and standard $L$--function. 
Hence, for $s \in \C$ with $\text{Re}(s) \gg 0$ the spinor $L$--function is denoted by $L(s,F)= \prod_p L_p(s,F)$
and the standard $L$--function by $L(s,F,st)= \prod_p L_p(s,F,st)$. These
are holomorphic functions. 
Let $F \in S_k$ be a primitive newform (normalized Hecke eigenform) with Fourier coefficients $a(n)$.
Then $L(s,F)=\sum_m a(m) \, m^{-s}$ and $L(s,F,st) = \frac{\zeta(2s-2k+2 )}{\zeta(s-k+1)} \sum_m a(m)^2 \, m^{-s}$.
We drop the index $n=1$ in the context of elliptic modular forms.

Let $g \in S_{2k}$ be a primitive newform with parameters $\{\alpha_p^{\pm}\}_p$
(determined by the eigenvalues $a(p)=(\alpha_p + \alpha_p^{-1}) p^{(2k-1)/2}$ 
and $\vert \alpha_p\vert =1$).
Then the $p$-part of the Hecke $L$--function $L(s,g)$ is given by 
$\left((1- \alpha_p p^{-s+k-(1/2)})(1- \alpha_p^{-1} p^{-s+k-(1/2)})\right)^{-1}$.
In \cite{Ik01}, for every positive integer $m \equiv k \pmod{2}$,
Ikeda explicitly constructed a Hecke eigenform $G = I_{k+m}^{2m}(g) \in
S_{k+m}^{2m}$, whose standard $L$--function is equal to 
\begin{equation}
\zeta(s) \prod_{j=1}^{2m} L(s+k+m-j,g).
\end{equation}
Recently Ikeda \cite{Ik06} introduced the concept of {\it Miyawaki
  lifts} and by techniques from representation theory  obtained the
following main result (Theorem 1.1 in his paper). 
Let $m,r,k \in \N$ and $n$ a non-negative integer with $m=n+r$ and $m \equiv k \pmod{2}$. Let $g \in S_{2k}$
be a primitive newform and $G\in S_{k+m}^{2m}$ be the Ikeda lift of $g$.
We fix the imbedding $j_1 \times j_2: \H_{2n+r} \times \H_r \longrightarrow \H_{2n+2r}$, given by
$(\tau,\tilde{\tau}) \mapsto 
\left( \begin{smallmatrix} \tau & 0 \\ 0 & \tilde{\tau} \end{smallmatrix} \right)$.
Then the pullback $G\circ (j_1\times j_2)$ can be viewed as a cusp form on $\H_r$ if we 
fix the variable on $\H_{2n+r}$.
Hence the function
\begin{equation}
\mathcal{G}_{g,h}(z):= \big\langle G \circ j_1(z,*), h \big\rangle,
\end{equation}
where $h \in S_{k+m}^r$ is well defined.
Let $g$ be a Hecke eigenform and $\mathcal{G}_{g,h}$ not identically zero.
Then $\mathcal{G}_{g,h}$ is a Hecke eigenform with standard $L$--function
\begin{equation}
L(s,h,st) \, \prod_{i=1}^{2n} L(s+k-n-i,g).
\end{equation}
\underline{Applications :}\\
Let $g_{20}\in S_{20}$ and $\Delta \in S_{12}$ be the unique primitive newforms.
By numerical calculations Ikeda has checked that
$\mathcal{G}_{g_{20},\Delta}$ is not identically zero. 
Hence this function has the same eigenvalues as the function $F_{20}$ introduced by
$Miyawaki$, since $\text{dim}S_{12}^3 = 1$. 
So, finally, the standard $L$--function of $F_{12}$ is equal to
\begin{equation}
L(s,\Delta,st) \, L(s+10,g_{20}) L(s+9,g_{20}).
\end{equation}
This proves Miyawaki's conjecture related to the standard $L$--function of $F_{12}$.

Ikeda's proof of Theorem 1.1 (see \cite{Ik06}) depends upon properties of
unramified principal series 
of $Sp_{2m}$ over $p$-adic fields. It is directly related to
the "even" Hecke algebra of the Hecke pair
$$
\big(   \Gamma_n, \text{M}_{2n}(\Z) \cap \text{Sp}_n(\Q)\big).
$$
Hence it does not cover the eigenvalues $\lambda_{F_{12}}(p)$ necessary to study the
spinor $L$--function of $F_{12}$.
To avoid this gap we choose a different approach to study the spinor
and standard $L$--functions 
of a Hecke eigenform, namely, we use Hayashida's description of Ikeda
lifts \cite{Ha07} in the frame work 
of Jacobi forms and Fourier-Jacobi expansions of Siegel modular forms.
It is in the spirit of Yamazaki's proof \cite{Ya86} of the Maass relations of Siegel Eisenstein series
and Eichler and Zagier's description of the Maass Spezialschar \cite{EZ85}.

Let $G \in S_{\kappa}^{2m}$ be a non-trivial cusp form of integer weight $\kappa$ and degree $2m$.
We will be mainly interested in the case $\kappa = k + n+1$ and $m=n+1$.
We further consider the Fourier Jacobi expansion
\begin{equation}
G \left( \begin{array}{cc} \tau & z^t \\ z & \tilde{\tau} \end{array}\right) = 
\sum_{l=1}^{\infty} 
\Phi_{\kappa,l}^G \left( \tau,z \right) \,\, \widetilde{q}^{\,\,l}
\end{equation}
where $\widetilde{q}:= e^{2 \pi \, i \, \tilde{\tau}}$.
It is well known that the Fourier Jacobi 
coefficients $\Phi_{\kappa,l}^G$ are Jacobi cusp forms 
$J_{\kappa,l}^{\text{cusp},2m-1}$ of weight $\kappa$, index $l$, and degree $2m-1$
on $\H_{2m-1}^J:= \H_{2m-1} \times \C^{2m-1}$.
Let $G= I^{2m}(g)$ be an Ikeda lift of the primitive newform $g\in S_{2k}$ with $\kappa = m +k$.
Then Hayashida \cite{Ha07} proved that there exists an operator 
\begin{equation}
D_{2m-1}\left( l, \{\alpha_p\}_p \right): J_{\kappa,1}^{\text{cusp},2m-1} 
\longrightarrow               J_{\kappa,l}^{\text{cusp},2m-1},
\end{equation}
where $\{\alpha_p\}_p$ are the parameters of $g$,
with the fundamental property
\begin{equation}
\Phi_{\kappa,l}^G = \Phi_{\kappa,1}^G\vert D_{2m-1}\left( l, \{\alpha_p\}_p \right).
\end{equation}
In the following we briefly recall the explicit definition of this operator.
We also would like to refer to the work of Katsurada and Kawamura \cite{KK06}
for a further application of the Hayashida's Fourier Jacobi model of Ikeda lifts.
Let $R$ be any subfield of $\R$. 
Let $$G^J_{2m-1}(R):= \left\{ \gamma= \left( \begin{array}{cccc} * & *&*&*\\ * & \alpha & * & * \\ 
* & *&*&*\\ 0 & 0 & 0 & \beta
\end{array}
\right) \in \text{G}^{+}\text{Sp}_{2m}(R)\big\vert \,\, \alpha,\beta \in R \text{ and positive}
\right\}.
$$
Then we set $\nu(\gamma):= \alpha/\beta$ and $\Gamma^J_{2m-1}:= \Gamma_{2m} \cap G^J_{2m-1}(\R)$.
With $\Phi \in J_{\kappa,l}^{\text{cusp},2m-1}$, put
$\widehat{\Phi}:= \Phi \,\,\widetilde{q}^l$. Define the Petersson slash operator for $\gamma \in 
G^J_{2m-1}(\R)$
via
$$\Phi\vert_{\kappa,l} \, \gamma := \left(\widetilde{q}\right)^{-l \nu(\gamma)}\, 
\widehat{\Phi}\vert_{\kappa}\gamma .$$
Yamazaki \cite{Ya86} introduced the two Hecke operators $T^J(p)$ and $T^J_{0,2m}(p^2)$ on the space 
$J_{\kappa,l}^{\text{cusp},2m-1}$ by the double cosets
\begin{equation}
\Gamma_{2m-1}^J \text{ diag}\left( 1_{2m-1},p ; p 1_{2m-1},1\right) \Gamma_{2m-1} \text{ and}
\end{equation}
\begin{equation}
\Gamma_{2m-1}^J \text{ diag}\left(  p 1_{2m-1},p^2 ; p 1_{2m-1},1\right) \Gamma_{2m-1} .
\end{equation}
More generally (we mainly follow Hayashida's normalization) let $M \in G_{2m-1}^J(\Z)$
and let the disjoint decomposition
\begin{equation*}
T^J := \Gamma_{2m-1}^J \, M \, \Gamma_{2m-1}^J = \bigcup_i \Gamma_{2m-1}^J M_i
\end{equation*}
be given. Then we define
\begin{equation}
\Phi \vert_{\kappa,l} T^J := \nu(M)^{m \kappa -(2m-1)m} \sum_i \Phi\vert_{\kappa,l} M_i.
\end{equation}
Now we recall the definition of the operators $D_{2m-1}\left( \{\alpha_p \}_p \right)$, where $\{\alpha_p\}_p$ are
the parameter of $g \in S_{2k}$ through the following formal Dirichlet series:
\begin{eqnarray*}
\quad \sum_{l=1} D_{2m-1} \left( l, \{\alpha_p \}_p \right) \,\, l^{-s}  & : = &\\
 & &  \!\!\!\! \!\!\!\!\!\!\!\! \!\!\!\!\!\!\!\! \!\!\!\!\!\!\!\! \!\!\!\! 
\!\!\!\! \!\!\!\!\!\!\!\! \!\!\!\!\!\!\!\! \!\!\!\!\!\!\!\! \!\!\!\!
\prod_p \big\{                                                                 
1 - G_p^{(m})(\alpha_p) \, T^J(p) \, p^{(m-1)(m+2)/2 -s} + T_{0,2m-1}^J(p^2) \, p^{2m(2m-1)-1-2s}
\big\}^{-1}.
\end{eqnarray*}
Here $G_p^{(1)}(\alpha_p) =1$ and
\begin{equation}
G_p^{(m)}(\alpha_p) = \prod_{ 1 \leq i \leq m-1}            \Big\{       
\left( 1 + \alpha_p p^{(1-2i)/2}\right)\,\, \left( 1 + \alpha_p^{-1} p^{(1-2i)/2}\right) \Big\}^{-1}
\end{equation}
otherwise.
Finally, we have
\begin{corollary}
Let $m \in \N$ and $p$ a prime. Then
\begin{equation}
D_{2m-1}( p, \{ \alpha_p\}_p) = G_p^{(m)}(\alpha_p) \,\, T^J(p) \, p^{(m-1)(m+2)/2}.
\end{equation}
\end{corollary}
In the cases $m=1$ and $m=2$ we have
\begin{equation}
D_1(p,\{ \alpha_p \}) = T^J (p) 
\text{ and } 
D_3(p, \{\alpha_p\}_p ) = p^2 \, G_p^{(2)}(\alpha_p) \,\, T^J(p).
\end{equation}
Moreover, let $b(n)$ be the Fourier coefficients of the primitive
newform $g \in S_{2k}$. Then we have 
$$G_p(\alpha_p) = \left\{p^k \left( b(p) + p^k +p^{k-1} \right)\right\}^{-1}.
$$

Let $\mathbb{W}$ be the Witt operator related to
the imbedding $\H_{2m-1} \times \H \longrightarrow \H_{2m}$, $\left(\tau, \tilde{\tau}\right) \mapsto
\left( \begin{smallmatrix}  \tau & 0 \\ 0 & \tilde{\tau} \end{smallmatrix}\right)$.
Then $\mathbb{W}\Phi_l^G$ is an element of $S_{k+m}^{2m-1}$ and
$\mathbb{W} G \in S_{k+m}^{2m-1} \otimes S_{k+m}$.
We consider the Witt map with respect to the Fourier Jacobi expansion of an 
Ikeda lift $G \in S_{\kappa}^{2m}$, where $m=n+1$ and $\kappa = k + n+1$.
$$
\mathbb{W}
\Big( \sum_{l=1}^{\infty} 
\Phi_{\kappa,l}^G \left( \tau,z \right) \,\, e^{2 \pi \, i \,l \tilde{\tau}}  \Big).
$$
Here $G$ is an Ikeda lift related to $g \in S_{2k}$ with Satake parameter $\alpha_p, \alpha_p^{-1}$ of
absolute value one.
Using Hayashida's discovery of a certain operator
$D_{2n-1}(m,\{\alpha_p\}_p)$ with the property 
$\Phi_{\kappa,l}^G = \Phi_1\vert D_{2m-1}(l,\{\alpha_p\}_p)$, our
approach is reduced to determination of
\begin{equation}
\sum_{l=1}^{\infty} 
\mathbb{W}
\Big( \Phi_{\kappa,1}^G\vert
D_{2m-1}(l,\{\alpha_p\}_p)\Big)
 \left( \tau,z \right) \,\, \widetilde{q}^{\,\, l}.
\end{equation}
We consider how the Witt map interacts with the operator 
$\vert_{\kappa} D_{2n-1}(m,\{\alpha_p\}_p)$. 
It is obvious that there exist symplectic Hecke operators $\mathcal{T}_l^{(2m-1)} \in \mathcal{H}^{2m-1}$ such that
\begin{eqnarray*}
\left( \mathbb{W} G  \right) (\tau,  \tilde{\tau}) & = & \sum_{l=1}^{\infty} 
\left( \mathbb{W}\Phi_{\kappa,l}^G \right) (\tau) \,\, \widetilde{q}^{\,\,l}\\
& = & \left(\mathbb{W} \Phi_{\kappa,1}^G\right) \vert_{\kappa}
\sum_{l=1}^{\infty} \mathcal{T}_l^{(2m-1)} 
 (\tau) \,\, \widetilde{q}^{\,\,l}.
\end{eqnarray*}

After a straightforward calculation we obtain the following Lemma.
\begin{lemma}
Let $G \in S_{\kappa}^{2m}$ be an Ikeda lift of $g\in S_{2k}$, where $\kappa=k+m$.
Then for every prime $p$:
\begin{equation}
\mathbb{W}\Phi_{\kappa,p}^G = 
p^{(m-1)(m+2)/2} \,\, p^{-(m-1)\kappa} \,\, G_p^{(m)}(\alpha_p)\,
\big( \mathbb{W}\Phi_{\kappa,1}^G\big) \vert_{\kappa} T^{(2m-1)}(p).
\end{equation}
\end{lemma}
\remark
For $m=1$ we have $\mathbb{W}\Phi_{\kappa,p}^G = 
\big( \mathbb{W}\Phi_{\kappa,1}^G \big) \vert_{\kappa} T^{}(p)$ and for $m=2$ we have
\begin{equation}
\mathbb{W}\Phi_{\kappa,p}^G = p^{2-\kappa}\,\, G_p^{(2)}(\alpha_p)\,\,
\big( \mathbb{W}\Phi_{\kappa,1}^G \big) \vert_{\kappa} T^{(3)}(p).
\end{equation}
In the more general case $m=3$ we have
\begin{equation}
\mathbb{W}\Phi_{\kappa,p}^G = p^{5-2\kappa}\,\, G_p^{(3)}(\alpha_p)\,\,
\big( \mathbb{W}\Phi_{\kappa,1}^G \big) \vert_{\kappa} T^{(5)}(p).
\end{equation}
%\newpage
%%
%%
%%
\begin{theorem}
Let $k >m$ and $G \in S_{k+m}^{2m}$ be an Ikeda lift attached to the primitive newform $g \in S_{2k}$.
For a Hecke eigenform $H \in S_{k+m}^{2m-1}$ and a primitive newform $h \in S_{k+m}$
we have: 
\begin{equation*}
\big\langle \mathbb{W}G, H \otimes h \big\rangle \neq 0
\end{equation*}
if and only if $\langle \mathbb{W}\Phi_{k+m,1}^G, H \rangle \neq 0$ and $\lambda_h(p)$ is equal to
$\lambda_H(p)$ times
\begin{equation}\label{form}
p^{(m-1)(m+2)/2} \,\, p^{-(m-1)(k+m)} \,\,
\prod_{ 1 \leq i \leq m-1}            \Big\{       
\left( 1 + \alpha_p p^{(1-2i)/2}\right)\,\, \left( 1 + \alpha_p^{-1} p^{(1-2i)/2}\right) \Big\}^{-1}
%\,\,
%\lambda_H(p) 
.
%%= \lambda_g(p).
\end{equation}
Here $\alpha_p, \alpha_p^{-1}$ are the parameters of $g$.
\end{theorem}

\begin{proof}
Assume $k>m$, so we can extend $H$ to an orthogonal Hecke eigenbasis 
$(H_i)_I$ of $S_{\kappa}^{2m-1}$, where $\kappa=k+m$ and
$H=H_1$.
Similiar we extend $h$ to a primitive newform basis $(h_j)_j$ of $S_{\kappa}$ with $h=h_1$.
Then we have
\begin{equation*}
\mathbb{W}G = \sum_{i,j} \alpha_{i,j}^G \, H_i \otimes h_j.
\end{equation*}
Hence, $\langle \mathbb{W}G,H_i \otimes h_j\rangle$ is equal to 
$\alpha_{i,j}^G \parallel H_i \parallel^2 \parallel h_j \parallel^2$.
On the other side, we can also employ the Fourier Jacobi model of
Ikeda lifts in the style of the Maass Spezialschar. Then we obtain for
$\langle \mathbb{W}G(\,\,,\tilde{\tau}),H_i\rangle$ the expression
\begin{equation}
\langle \mathbb{W}\Phi_{\kappa,1}^G, H_i \rangle \,\, 
\sum_{l=1}^{\infty} \lambda_{H_i}( \mathcal{T}_l^{2m-1}) \,\, \widetilde{q}^l.
\end{equation}
Here we used several results related to the Fourier Jacobi expansion
of Ikeda lifts and the fact that the operators $\mathcal{T}_l^{2m-1}$
had been choosen and normalized in such a way that they are
self-adjoint with respect to the Petersson scalar product. 
For $l=1$ we have the identity.
Let $\langle \mathbb{W}\Phi_{\kappa,1}^G, H_i \rangle \neq 0$. Then it
follows from \cite{Ik06} that 
$\sum_{l=1}^{\infty} \lambda_{H_i}( \mathcal{T}_l^{2m-1}) \,\, \widetilde{q}^l$
is a primitive newform, since multiplicity-one for $\text{SL}_2$ is available \cite{Ra00}.
Hence, it is equal to one of the newforms among the basis $(h_j)_j$.
Since the eigenvalues $\lambda_{H_i}( \mathcal{T}_p^{2m-1})$ already determine the newform,
the theorem is proven.
\end{proof}
\begin{corollary}
Let $G \in S_{\kappa}^{2m}$ be an Ikeda lift with non-trivial Witt image $\mathbb{W}G$.
Then there exists a Hecke eigenform $H \in S_{\kappa}^{2m-1}$ and a primitive newform $g \in S_{\kappa}$
such that $a_h(p)$ is equal to $\lambda_H(p)$ times the expression (\ref{form}).
If we normalize $H$ in the case $m=1$, then $H=h$.
\end{corollary}
Let $G \in S_{\kappa}^4$ be an Ikeda lift ($\kappa=k+2$) associated to the primitive newform $g \in S_{2k}$
with non-trivial Witt image.
Then there exists a Hecke eigenform $H \in S_{\kappa}^3$ and a
primitive newform $h\in S_{\kappa}$ such that
\begin{equation}
\lambda_H(p) = \lambda_h(p) \Big( \lambda_g(p) + p^k + p^{k-1}\Big).
\end{equation}
{\underline{Applications :}}\\
Let $G \in S_{12}^4$ be the Ikeda lift associated to $g_{20} \in S_{20}$.
Then numerical observations show that the Witt image of $F$ is
non-trivial (see for example \cite{Ik06}). 
Since $\text{dim}S_{12}^3 = \text{dim}S_{12} =1$ we obtain
\begin{equation}
\lambda_{F_{12}}(p) = a_{\Delta}(p) \Big( a_{g_{20}}(p) + p^{10} + p^9 \Big).
\end{equation}
%%\newpage
\section{The spinor $L$--function of Miyawaki lifts of degree $3$}
Miyawaki formulated his conjectures very explicitly.
Hence, we restrict ourself in this section exclusively to
Siegel cusp forms of degree $3$.
We do not consider Eisenstein series, since their $L$--functions
are related to modular forms of lower degree via the Siegel $\Phi$-operator.

Let $F \in S_{\kappa}^3$ be a Miyawaki lift associated to an Ikeda lift of degree $4$.
We prove that the spinor $L$--function of $F$ satisfies
the Miyawaki property (P) given in \cite{Mi92}, page 326.
This can be stated in the following way:
There exist two primitive newforms $f \in S_{k_1}$ and $g\in S_{k_2}$ with
$k_1= \kappa$ and $k_2 = 2 \kappa-4$ such that the
spinor $L$--function $L(s,F)$ is given by
\begin{equation}\label{deg}
L(s,F) = L(s - \kappa +2,f) \,\, L(s -\kappa+3,f) \,\, L(s,f \otimes g).
\end{equation}
Here $L(s,f)$ is the Hecke $L$--function of $f$ and $L(s,f \otimes g)$
the Rankin $L$--function of $f$ and $g$.
For the readers convenience we give a precise definition of these $L$--functions.
Let $h \in S_k$ be a primitive newform with local parameters $\{\alpha_p^h\}_p$. Then
\begin{equation}
L(s,h) := \sum_{m=1}^{\infty} \,\, a_h(m) \, m^{-s} = \prod_p \text{det} 
\Big( 1_2 - p^{(k-1)/2} \left( \begin{array}{cc} \alpha_p & 0 \\ 0 & \alpha_p^{-1}
\end{array}
\right)\Big)^{-1}
\end{equation}
is a the Hecke $L$--function of $h$.
Further, let $f \in S_{k_1}$ and $g \in S_{k_2}$ be primitive
newforms. Then the Rankin convolution of 
$f$ and $g$ is given by
\begin{equation}
L(s,f \otimes g) := 
\prod_p \text{det} 
\Big( 1_4 - p^{(k_1+k_2-2)/2} 
\left( \begin{array}{cc} \alpha_p^f & 0 \\ 0 & \left(\alpha_p^f\right)^{-1}
\end{array}\right)
\otimes
\left( \begin{array}{cc} \alpha_p^g & 0 \\ 0 & \left(\alpha_p^g\right)^{-1}
\end{array}\right)
\Big)^{-1}.
\end{equation}
More generally, let $\mu_0,\mu_1, \ldots, \mu_n$ be the $p$-Satake parameter of 
the Hecke eigenform $F \in S_k^n$. These are complex numbers, which are unique
up to the action of the Weyl group of the symplectic group and
$\mu_0^2 \mu_1 \cdot \ldots \cdot \mu_n = p^{nk -n(n+1)/2}$.
This is compatible with our normalization of the Hecke operators.
Because of the Satake isomorphism
we can also define the spinor and the standard $L$--function
via these parameters. The local $L$--factors are
\begin{eqnarray}
L_p(s,F) & :=& \left( 1 - \mu_0 p^{-s}\right)^{-1} \prod_{r=1}^n \prod_{i_1 < \ldots <i_r} 
\left( 1 - \mu_0 \mu_{i_1} \cdot \ldots \cdot \mu_{i_r} p^{-s}\right)^{-1},\\
L_p(s,F,st) & := & \left( 1-p^{-s}\right)^{-1} \prod_{i=1}^n 
\left( 1 - \mu_i p^{-s}\right)^{-1}
\left( 1 - \mu_i^{-1} p^{-s}\right)^{-1}.
\end{eqnarray}
Then the spinor $L$--function $L(s,F)$ is equal to $\prod_p L_p(s,F)$ and the
standard $L$--function is defined via
\begin{equation}
L(s,F,st) := \prod_p L_p(s,F,st).
\end{equation}
The spinor $L$--function can also be defined by formal power series related to
the Hecke operators $T(p^d)^{(n)}$ in the following way:
\begin{equation}
\sum_{d=0}^{\infty} \, T(p^d)^{(n)} \, X^d = \frac{P_p(X)}{Q_p(X)},
\end{equation}
where $P_p(X)$ and $Q_p(X)$ are polynomials of degree $2^n-2$ and $2^n$
with coefficients in the Hecke algebra.
The denominator polynomial $Q_p(X)$ is directly related to
the local $L$--factor of the spinor $L$--function.
In the case $n=1$ we have
\begin{equation}
Q_p(X) = 1 -T(p)X + p T_1(p^2) X^2.
\end{equation}
If we replace the operators by the eigenvalues 
we obtain the polynomial $Q_{p,F}(X)$. Here we drop the index $(n)$ to simplify notation.
Hence
$L_p(s,F) = Q_{p,F}(p^{-s})^{-1}$. For $n=2$ we have
\begin{equation*}
1-T(p)X + \left( T(p)^2 + p(p^2 +1)T_2(p^2)
\right)X^2 - p^3 T(p) T_2(p^2) X^3 + p^6 (T_2(p^2))^2 X^4.
\end{equation*}
The case $n=3$ is complicated. It had been first given
by Andrianov \cite{An67}, in connection with a proof of the conjecture of Shimura
related to the symplectic group of genus $3$.
We have $Q_p(X) = \sum_{m=0}^8 (-1)^m \, c(m) \, X^m$ with
\begin{eqnarray}
c(0) & = &    1 ,         \label{c0}                 \\
c(1) & = &    T(p) ,            \label{c1}               \\
c(2) & = &      p\Big(   T_1(p^2) + (p^2 +1) T_2(p^2) + (p^2+1)^2 T_3(p^2)\Big) ,           \label{c2}               \\
c(3) & = &     p^3 T(p) \Big( T_2(p^2) + T_3(p^2) \Big)  ,       \label{c3}                   \\
c(4) & = &     p^6 \Big(
T(p)^2 T_3(p^2) + T_2(p^2)^2 - 2p T_1(p^2) T_3(p^2) 
\label{c4} \\
& &- 2 (p-1) T_2(p^2) T_3(p^2)-(p^6 + 2p^5 + 2p^3 + 2p-1) T_3(p^2)^2\Big),
 \nonumber                          \\
c(5) & = &  p^6 T_3(p^2) c(3)  ,        \label{c5}                     \\
c(6) & = &  p^{12} T_3(p^2)^2 c(2),       \label{c6}                        \\
c(7) & = &   p^{18} T_3(p^2)^3 c(1),    \label{c7}                          \\
c(8) & = &   p^{24} T_3(p^2)^4.    \label{c8}                          
\end{eqnarray}
We have listed these coefficients explicitly, and double-checked them,
since we need the exact values to prove the Miyawaki conjecture. 

The degeneration of the spinor $L$--function $L(s,F)$ given in (\ref{deg}) of
the Hecke eigenform $F \in S_{\kappa}^3$, where $\kappa = k+2$, is
equivalent to the following formulas of the eigenvalues:
\begin{eqnarray*}
\lambda_F(p) & = & a_f(p) \Big( a_g(p) + p^{k_2/2} + p^{k_2/2-1} \Big),\\
\lambda_F(T_3(p^2)) &= &  p^{3\kappa-12},\\
\lambda_F(T_2(p^2)) & = &
a_f(p)^2 p^{k_2-4} + a_g(p) p^{k_1 + k_2/2-5} (p+1) - p^{3\kappa-12}(p^3 +1),\\
\lambda_F(T_1(p^2)) & = &
a_f(p)^2 a_g(p) p^{k_2/2 -2 } (p +1) + a_f(p)^2 p^{k_2 -4}(p^2-1) 
\\
& &
+ a_g(p)^2 p^{k_1-2} - a_g(p) p^{k_1+k_2/2 -5} (p^2+1)(p+1) \\
& &+ p^{3\kappa-10}(p^3+1)(p-1).
\end{eqnarray*}
We have already proven the formula for $\lambda_F(p)$.
For this we employed the Fourier-Jacobi model of Ikeda lifts
introduced by Hayashida, in the style of the Maass lifts in the setting
of Saito-Kurokawa lifts. Moreover we
discovered Hecke duality properties of the Witt operator on Ikeda lifts of degree $4$
and the involved pullback components.
The formula for $\lambda_F(T_3(p^2))$ is obvious and does not depend on
the Miyawaki property.
The operator $T_3(p^2)$ has only one left coset $\Gamma_3 \, p 1_6$.
Hence $F\vert_{\kappa} T_3(p^2)$ is equal to $p^{3 \kappa -12}F$.

Next we want to prove the formula for the eigenvalue $\lambda_F(T_2(p^2)$.
This may be possible by studying the formal power series
%\begin{equation*}
$\sum_{l=1}^{\infty} \mathcal{T}_l^{(3)} \,\, X^l
%\end{equation*}
$
in more detail.
We choose a different way which works for Miyawaki lifts 
for degrees larger than $3$, too.

Recall the following notation. Let $F \in S_{\kappa}^3$ be the
Miyawaki lift of the Ikeda lift $G \in S_{\kappa}^4$ 
and the primitive newform $f \in S_{\kappa}$ be given. Here $\kappa = k+2$ and
$G$ is the lift of the primitive newform $g \in S_{2k}$.
Let $\{\alpha_p^g\}_p$ be the parameters of $g$ and $\left( \beta_{p,0}^f, \beta_{1,p}^f\right)_p$
be the Satake parameters of $f$. We drop the index $p$ to simplify notation.
Further, the Satake parameters of $G$ related to the standard zeta function
had been determined by Ikeda \cite{Ik06} using representation theory.
Let $\mu_0^G,\mu_1^G,\mu_2^G,\mu_3^G$ be the Satake parameter of $G$.
Then we can choose uniquely (up to the action of the symplectic Weyl group):
\begin{eqnarray*}
\mu_1^G &=& \beta_1^f,\\
\mu_2^G & = & \alpha^g \, p^{1/2},\\
\mu_3^G & = & \left(\alpha^g\right)^{-1} \, p^{1/2}.
\end{eqnarray*}
These equations determine $\mu_0^G$ for every prime $p$ up to sign.
\begin{equation}
\left(\mu_0^G\right)^2 = p^{3 \kappa - 7}/\beta_1^f.
\end{equation}
Next we claim that we can determine the Miyawaki formula for $T_2(p^2)$ from
formula (\ref{c1}) and (\ref{c3}).
Let $S(a,b,c):= 1+a +b +c + ab + ac + bc + abc$, then
\begin{equation}
\mu_0^G \, S\Big(\beta_1^f, \alpha^g p^{1/2}, \left(\alpha^g\right)^{-1} p^{1/2} \Big) = \lambda_F(p).
\end{equation}
Further, let the polynomial $T(a,b,c)$ be given by:
$$
\begin{array}{l}
ab(c + ab + ac + bc + abc) + 
ac(ab + ac + bc + abc) + aab(ac + bc + abc)\\
+ 
aac(bc + abc) + abcabc+
bc(ab + ac + bc + abc)+bab(ac + bc + abc)\\+ 
b ac(bc + abc)+ b bcabc +
cab(ac + bc + abc) + cac(bc + abc)
\\+ cbcabc +
abac(bc + abc)+abbcabc+
acbcabc+
a(b+c+ab+ac+bc+abc)\\
+b(c+ab+ac+bc+abc)+
c(ab+ac+bc+abc)+
ab(ac+bc+abc)\\+
ac(bc+abc)+
bcabc.
\end{array}
$$
Then we obtain from (\ref{c4}) the equation
\begin{eqnarray}
\left(\mu_0^G\right)^2 \, 
T\Big(\beta_1^f, \alpha^g p^{1/2}, \left(\alpha^g\right)^{-1} p^{1/2} \Big) & & \nonumber\\
&    & 
\!\!\!\!\!\!\!\!\!\!\!\!\!\!\!\!\!\!\!\!\!\!\!\!\!\!\!\!\!\!\!\!\!\!\!\!
\!\!\!\!\!\!\!\!\!\!\!\!\!\!\!\!\!\!\!\!\!\!\!\!\!\!\!\!\!\!\!\!\!\!\!\!
=
p^3 \, S\Big(\beta_1^f, \alpha^g p^{1/2}, \left(\alpha^g\right)^{-1} p^{1/2} \Big) 
\Big( \lambda_F(T_2(p^2)) + p^{3 \kappa -12}\Big).
\end{eqnarray}
Let $\lambda_F(p) \neq 0$ then we can conclude
\begin{equation}
\lambda_F(T_2(p^2)) = p^{3 \kappa -12} \Bigg( \frac{p^2}{\beta_1^f} 
\frac{T\big(\beta_1^f, \alpha^g p^{1/2}, \left(\alpha^g\right)^{-1} p^{1/2} \big)
}{
S\big(\beta_1^f, \alpha^g p^{1/2}, \left(\alpha^g\right)^{-1} p^{1/2} \big) 
} 
-1 \Bigg).
\end{equation}
After a straightforward calculation we obtain for the right side of
this equation
\begin{eqnarray*}
p^{3 \kappa-9} \left( \beta_1^f + (\beta_1^f)^{-1} +2 \right) 
+
\left(\alpha^g + (\alpha^g)^{-1}\right)
p^{(2 \kappa -5)/2} \, p^{2 \kappa -7}(p+1) - p^{3\kappa-12}(p^3+1).
\end{eqnarray*}
On the other side the Fourier coefficients $a_f(p),a_g(p)$ are related to
the above parameters. We have
\begin{equation}
a_f(p)^2 = p^{\kappa-1} \left( \beta_1^f + (\beta_1^f)^{-1} +2 \right) \text{ and }
a_g(p) = p^{(2 \kappa-5)/2} \left(\alpha^g + (\alpha^g)^{-1}\right).
\end{equation}
Hence the conjecture of Miyawaki has been proven for $\lambda_F(T_2(p^2))$ if the
the eigenvalue $\lambda_F(p) \neq 0$.
Similiar we obtain the expected formula for $\lambda_F(T_1(p^2))$ if we evaluate the equation
(\ref{c2}) in two ways and equate them.
One way is via the description by Satake parameters and the other one
is by plugging in 
the formula already obtained for $T_2(p^2)$.

So finally it remains to treat the degenerate case $\lambda_F(p)=0$
which perhaps does not occur, but which we cannot omit, since
the generalized Lehmer conjecture has not been proven yet.
We have 
$\lambda_F(p) = a_f(p) \left( a_g(p) + p^{k_2/2} + p^{k_2/2-1}\right)$.
But since the primitive newform $g \in S_{k_2}$ 
has the property $\vert a_g(p)\vert \leq 2 p^{\frac{k_2 -1}{2}}$
we have $\lambda_F(p) = 0$ if and only if $a_f(p)=0$.
Assume that $a_f(p)=0$. Then the Satake parameters $\left( \mu_0^F, \mu_1^F,\mu_2^F, \mu_3^F\right)$
of $F$ are given by
\begin{equation}
\Big(
\varepsilon (\beta_1^f)^{-1/2} p^{(3 \kappa-7)/2}, 
\beta_1^f, \alpha^g p^{1/2}, (\alpha^g)^{-1} p^{1/2} \Big),
\end{equation}
with $\beta_1^f = -1$ and $\varepsilon = \pm 1$.
Since the Satake parameters are only unique up to the action of
the related symplectic Weyl group, 
we can choose $\varepsilon=1$, because $\beta_1^f$ is degenerate and has
the value $(-1)$. 
The $p$-local factor of $L(s,F)$ in this special case can be
carried out directly via the Satake parameters.
This leads to the predicted formulas 
for the spinor L-function and of the eigenvalues of $F$.
\subsection{Final remarks}
In the paper \cite{Mu02}, Murakawa gives two theorems related to the spinor L-function
of Ikeda lifts (Theorem 3.1) and Miyawaki lifts (Theorem 5.1). Since our approach is different
it is maybe worthwile to make some comments.

The Satake parameters of the related standard L-functions are given by Ikeda \cite{Ik01},\cite{Ik06}.
To talk about spinor L-function one has first to show that Ikeda lifts and 
Miyawaki lifts are Hecke eigenform for the full Hecke algebra. This is not proven in Ikeda's paper \cite{Ik01}, since
Ikeda had been only interested in the Hecke algebra related to $\text{S}p_{2n}{(\Q)}$ 
and not $\text{G}^{+}\text{Sp}_{2n}(\Q)$. 
Further one has to determine the sign of the Satake parameter $\mu_0$ in the case of Ikeda and Miyawaki lifts
(for Miyawaki lifts one also has to care about the non-vanishing).

Murakawa first gives a proof of Theorem 3.1. He assumes that $\mu_0$ exists.
Since Murakwa's paper is an extraction of his master thesis in Japanese we may assume that 
the related proof is given there. Then he calculates $\mu_0$ by a method indicated 
by Ikeda in the case of the standard L-function. This works since the related Satake parameter can be obtained
via equating certain equations in which the involved terms are none-zero.

After stating Theorem 5.1 he indicated that the proof can be given in the same way as the proof of Theorem 3.1,
if one exchanges the Siegel type Eisenstein series by a Klingen type Eisenstein series. Again one
first would have to show that the Miyawaki lift is a Hecke eigenform for the full Hecke algebra
(in our proof we had to use the multiplicity one theorem of $\text{SL}_2$, hence one may could give a new
prove of this by applying Murakawa's result, if fully available).
Moreover we found out that only if Lehmers conjecture is true one could tranfer Lemma 4.1 \cite{Mu02}, since
in contrast to Ikeda lifts the eigenvalues $\lambda(p)$ of Miyawaki lifts could be zero.
In the case of $F_{12}$ we actually proved that all $\lambda(p)$ are non-zero if and only if the Lehmer conjecture is true.
That means even when the Satake paramters $\mu_1,\mu_2,\mu_3$ of the Miyawaki lift $F_{12}$ are fixed it may happen
that $\mu_0$ is not unique.

\end{document}